\documentclass{article}
\usepackage{amssymb, amsmath, amsthm}
\usepackage[arrow,matrix]{xy}
\usepackage[dvips]{epsfig}
\usepackage[active]{srcltx}
\newtheorem{theorem}{Theorem}
\newtheorem{proposition}{Proposition}
\newtheorem{corollary}{Corollary}
\newtheorem{lemma}{Lemma}

\newcommand\R{\mathbb R}
\newcommand\Z{\mathbb Z}
\newcommand\C{\mathbb C}
\def\sgo{singular Green operator}
\def\cstar{$C^*$-algebra}

\def\p0#1{\overline{\Psi^0}(#1)}
\newcommand\adg{{}^{ad}\kern-1pt \mathcal{G}}
\def\g0{G^{(0)}}
\newcommand\ad[1]{{}^{ad}\kern-1pt #1}

\title{Boutet de Monvel's Calculus and Groupoids I}
\author{J. Aastrup, S. T. Melo\thanks{S.T. Melo and B. Monthubert were partly funded by a
  cooperation agreement CAPES-COFECUB.}, B. Monthubert \&\ E. Schrohe}
\date{}
 
\begin{document}
\maketitle

\begin{abstract}
Can Boutet de Monvel's algebra on a compact manifold with
boundary  be obtained as the algebra $\Psi^0(G)$ 
of pseudodifferential operators on some Lie groupoid $G$? 
If it could, the kernel ${\mathcal G}$ of the principal symbol homomorphism 
would be isomorphic to the groupoid \mbox{$C^*$-algebra}  $C^*(G)$. 
While the answer to the above question remains open,  we exhibit in this paper a groupoid 
$G$ such that $C^*(G)$ possesses an ideal ${\mathcal I}$ isomorphic to ${\mathcal G}$.
In fact, we prove first that ${\mathcal G}\simeq\Psi\otimes{\mathcal K}$ with the $C^*$-algebra 
$\Psi$ generated by the zero order pseudodifferential operators on the boundary and the 
algebra $\mathcal K$ of compact operators. 
As both $\Psi\otimes \mathcal K$ and $\mathcal I$ are extensions of 
$C(S^*Y)\otimes {\mathcal{K}}$  by ${\mathcal{K}}$ ($S^*Y$ is
the co-sphere bundle over the boundary) we infer from 
a theorem by Voiculescu that both are isomorphic.
\end{abstract}

\section*{Introduction} 
Boutet de Monvel's calculus \cite{B1-1, B1-2,bdm,grubb,rempel-schulze,schrohe} is a pseudodifferential calculus on manifolds with boundary.
It includes the classical differential boundary value problems as well as the parametrices to elliptic elements. 
Many operator-algebraic aspects of this algebra (spectral invariance, noncommutative
residues and traces, composition sequence, $K$-theory) have been studied
recently \cite{fedosov-golse-leichtnam-schrohe,grubb-schrohe,melo-nest-schrohe1,melo-schick-schrohe,nest-schrohe,schrohe99}. The problem of identifying this
algebra as the pseudodifferential algebra (or as an ideal of one) of a Lie groupoid may be the key to
an effective application of the methods of noncommutative geometry. 
If that is acomplished, one could then seek for extensions of the calculus,
and for a better understanding of its index theory \cite{bdm,fedosov,rempel-schulze}. Basic definitions and certain facts about Boutet de Monvel's
algebra are recalled in Section~\ref{bdmc}. 

The groupoid approach to pseudodifferential calculus was
developed in noncommutative geometry, following the seminal work of
A. Connes for foliations \cite{connes}. The guiding principle in that approach 
is that the central object in global analysis is a groupoid. 
To study a particular situation, for a class of singular manifolds for
instance, one needs to define a groupoid adapted to the context and
then use the general pseudodifferential tools for groupoids, as 
developed in \cite{bm-fp, nwx, bm-jfa, lmn, lmn2}. 
This has been applied to the context of manifolds with corners, with
the goal of studying Melrose's $b$-calculus (see \cite{bm-fp, nwx,
    bm-jfa}). Groupoids were defined whose pseudodifferential calculi
  recover the $b$-calculus and the cusp-calculus.
Basic definitions and certain facts about pseudodifferential calculus on
groupoids are recalled in Section~\ref{psi}. 

The starting point of this work is the following result (see
\cite{melo-nest-schrohe1}): 
The kernel of the 
principal symbol map for Boutet de Monvel's 
calculus is equal to the norm closure $\mathcal G$ 
of the ideal of singular Green operators. 
Since in the pseudodifferential calculus on a groupoid, the
$C^*$-algebra of the groupoid is the kernel of the principal symbol
map, this gives a hint that finding a groupoid whose $C^*$-algebra is
$\mathcal{G}$ could give some insight about the relationship between
the Boutet de Monvel algebra and groupoid pseudodifferential algebras.

Besides, $\mathcal G$
fits into a 
short exact sequence (see  \cite[Section 7]{melo-schick-schrohe}):
$$0 \to \mathcal{K} \to \mathcal{G} \to C(S^*Y) \otimes \mathcal{K}
\to 0,$$
which is similar to that for pseudodifferential operators
on smooth manifolds:
$$0 \to \mathcal{K} \to \Psi \to C(S^*Y)
\to 0.$$
 In Section~\ref{sgoi}, we show
that ${\mathcal G}$ is actually Morita-equivalent to the 
norm-closure $\Psi$ of the algebra of pseudodifferential operators 
on the boundary.
Since ${\mathcal G}$ is a stable \cstar, 
it is thus isomorphic
to $\Psi\otimes \mathcal{K}$.

On the other hand, we define in Section \ref{gr-sgo} a groupoid
whose $C^*$-algebra contains an ideal $\mathcal{I}$ which fits
in an extension analoguous as that of $\Psi\otimes \mathcal{K}$. By
showing in Section \ref{main} that the $KK$-theory elements induced by these extensions
coincide, we infer from a theorem by Voiculescu 
that ${\mathcal G}$ and $\mathcal{I}$ are isomorphic. 

%
%
%
%

\section{Boutet de Monvel's Calculus}\label{bdmc}
Let $X$ denote a compact manifold of dimension $n$ with boundary $Y$ and
interior $\dot X$. Given a pseudodifferential operator $P$, defined on
an open neighborhood $\tilde X$ of $X$, and $u\in C^\infty(X)$, one
defines $P_+u$ as equal to the restriction to $\dot X$ of $PEu$, where 
$Eu$ is the extension of $u$ to $\tilde X$ which vanishes outside
$X$. In general, singularities may develop at the boundary, and one
gets only a mapping $P_+:C^\infty(X)\to C^\infty(\dot X)$. One says 
that $P$ has the {\em transmission property} if the image of the truncated operator $P_+$ is 
contained in $C^\infty(X)$ (a subspace of $C^\infty(\dot X)$). This
was defined by Boutet de Monvel in \cite{B1-1, B1-2}, 
where he proved that the transmission property for a classical 
(polyhomogoneous) pseudodifferential operator is 
equivalent to certain symmetry conditions 
for the homogeneous components of the asymptotic expansion of its symbol at 
the boundary.
Later \cite{bdm}, he constructed an algebra whose elements are operators of the form
\begin{equation}
\label{bdmo}
\left(
\begin{array}{cc}
P_+\,+\,G&K\\\\T&S
\end{array}
\right)
\,:\,
\begin{array}{c}
C^\infty(X)\\\oplus\\C^\infty(Y)
\end{array}
\longrightarrow
\begin{array}{c}
C^\infty(X)\\\oplus\\C^\infty(Y)
\end{array},
\end{equation}
where $P$ is a pseudodifferential operator on $X$ satisfying a condition that ensures the transmission
property, $S$ is a pseudodifferential operator on $Y$, while $G$, $K$
and $T$ belong to classes of operators he then defined and named,
respectively, singular Green, Poisson and trace operators. We call an operator as in (\ref{bdmo}) a 
Boutet de Monvel operator. For detailed expositions of his calculus, we refer to \cite{grubb,rempel-schulze}. 

A Boutet de Monvel operator has an {\em order},
roughly the usual order of pseudodifferential operators. The entries $T$ and $G$
have, moreover, an integer {\em class} assigned to them. The class of a trace operator is related to the order of the derivatives that
appear in the boundary-value conditions it prescribes. One must
assign a class also to singular Green operators due to the fact that
the composition $KT$ is a singular Green operator. There exist isomorphisms between suitable
Sobolev spaces such that the composition of a given operator of
arbitrary order and class with one of them has order and class
zero. For index theory purposes it is therefore sufficient to consider
operators of order and class zero. 
These form an adjoint invariant
subalgebra of the algebra ${\mathcal L}({\mathcal H})$ of all bounded operators on a suitable Hilbert
space ${\mathcal H}$. 
Adopting the definition of order in  \cite{rempel-schulze,schrohe} 
for $K$ and $T$, we here choose  ${\mathcal H}=L^2(X)\oplus L^2(Y)$. 
If, as does Grubb \cite{grubb}, one keeps the original 
definition (which makes more sense if one needs general $L^p$ estimates) 
then one must take a Sobolev space of order $-1/2$ over the boundary.

Boutet de Monvel operators can also be defined as mappings between smooth sections
of vector bundles. If $E$ is a bundle of positive rank over $X$, and $F$ is an arbitrary bundle 
over $Y$, then the algebra of all Boutet de Monvel operators of order and class zero acting between sections of $E$ and
$F$ is Morita equivalent \cite[Section 1.5]{melo-nest-schrohe1} to the algebra obtained by taking a 
rank-one trivial bundle over $X$ and the zero-bundle over $Y$. 
This partly justifies, again if one is interested in index theory, 
to consider only the operators appearing in the 
upper-left corner of the matrix in (\ref{bdmo}) and to assume, as we did at the beginning, 
that the bundle over $X$ is $X\times {\mathbb C}$.

The problem of computing the Fredholm index of a Boutet de Monvel 
operator acting between sections of different bundles over each side can be reduced to the case of equal bundles on both 
sides by a device developed by Boutet de Monvel \cite{B1-1, B1-2}, recalled in \cite[Section 1.1]{melo-schick-schrohe}.

Let us now explain what a singular Green operator $G$ is, in the case of order and class zero and of a 
rank-one trivial bundle over $X$. Its distribution kernel is smooth outside
the boundary diagonal; i.e, if $\varphi\in C_c^\infty(\dot X)$, and if we
denote by $M_\varphi$ the operator of mulitiplication by $\varphi$, then $GM_\varphi$ and $M_\varphi G$ 
are integral operators with smooth kernels. The push-forward of $G$ by a boundary chart is
an operator-valued-symbol pseudodifferential operator on the variables tangential
to the boundary, as we describe below. 
It is perhaps worth stressing, however, that it is in general not a pseudodifferential 
operator on all variables, because of its particular way of acting on the normal variable. 

Given $u\in C_c^\infty({\mathbb R}^n_+)$, ${\mathbb R}_+^n
=\{(x^\prime,x_n)\in{\mathbb R}^{n-1}\times{\mathbb R};\ 
x_n\geq 0\}$, let $\hat u$ denote the vector-valued Fourier transform of $u$ with
respect to the $n-1$ first variables, 
\begin{equation}
\label{ft}
\hat u(\xi^\prime)=
\int e^{ix^\prime\cdot\xi^\prime}u(x^\prime,\cdot)dx^\prime \in C_ c^\infty({\mathbb R}_+).
\end{equation}
In local coordinates for which the boundary corresponds to $x_n=0$ and the interior to $x_n>0$, $G$ is given by
\begin{equation}
\label{defg} 
Gu(x^\prime,\cdot)=(2\pi)^{1-n}\int
e^{ix^\prime\cdot\xi^\prime}g(x^\prime,\xi^\prime,D_n)\hat u(\xi^\prime)\,d\xi^\prime
\end{equation}
The integrals in (\ref{ft}) or in (\ref{defg}) should be regarded, for fixed $\xi^\prime$ or $x^\prime$, respectively, 
as $L^2({\mathbb R}_+)$-valued integrals. For each $(x^\prime,\xi^\prime)$, $g(x^\prime,\xi^\prime,D_n)$ in (\ref{defg})
is an integral operator with kernel $\tilde g(x^\prime,\cdot,\cdot,\xi^\prime)$ equal to the restriction to
${\mathbb R}_+\times{\mathbb R}_+$ of a function belonging to the
Schwartz space of rapidly decreasing functions on ${\mathbb R}^2$. The
function $\tilde g(x^\prime,x_n,y_n,\xi^\prime)$ (called by Grubb the
{\em symbol-kernel} of $G$) is smooth and satisfies the estimates \cite[(1.2.38)]{grubb}. 
This is invariantly defined \cite[Theorem 2.4.11]{grubb} with respect to coordinate changes 
that preserve the set $\{x_n=0\}$.

We denote by ${\mathcal A}_0$ the set of all polyhomogeneous operators $P_++G$ of order and class zero on 
$X$, and by ${\mathcal G}_0$ its subset of all singular Green operators. It follows from the rules of Boutet 
de Monvel's calculus that ${\mathcal A}_0$ is an algebra and  
that ${\mathcal G}_0$ is an ideal in ${\mathcal A}_0$. 

In the sequel, we shall restrict ourselves to coordinate changes which
preserve the variable $x_n$, i.e., we choose a normal coordinate.
Then  two *-homomorphisms are defined on ${\mathcal A}_0$, 
the principal symbol and the boundary principal symbol:
\[
\sigma:{\mathcal A}_0\to C^\infty(S^*X)\ \ \mbox{and}\ \ \gamma:{\mathcal A}_0\to C^\infty(S^*Y,{\mathcal L}(L^2({\mathbb R}_+))). 
\]
The principal symbol of a given $P_++G$ is, by definition, the usual principal symbol of $P$
\[
\sigma(P_++G)=p_0,
\]
where $p_0$ is the leading term in the aymptotic expansion of the symbol of $P$. 

At a point $(x',\xi')$ in $S^*Y$, the boundary principal symbol of  $P_+$
is defined to be the truncated Fourier multiplier
\[
\gamma_{P_{+}}(x^\prime,\xi^\prime)=p_0(x^\prime,0,\xi^\prime,D_n)_+
\]
of symbol $\xi_n\mapsto p_0(x^\prime,0,\xi^\prime,\xi_n)$. 
The boundary principal symbol of $G\in\mathcal G_0$ is the integral operator
\begin{equation}\label{bpsg}
\gamma_G(x^\prime,\xi^\prime)=g_0(x^\prime,\xi^\prime,D_n)
\end{equation}
with the rapidly decreasing kernel 
$\tilde g_0(x^\prime,\cdot,\cdot,\xi^\prime)$, where $\tilde g_0$ denotes the leading term
in the asymptotic expansion of $\tilde g$, cf. \cite[(1.2.39)]{grubb}. 
Then  $\gamma$ maps ${\mathcal G}_0$ into 
$C^\infty(S^*Y,{\mathcal{K}}_{{\mathbb R}_{+}})$, 
with the ideal ${\mathcal{K}}_{{\mathbb R}_{+}}$  of compact 
operators on $L^2({\mathbb R}_+)$.

Let ${\mathcal A}$ and ${\mathcal G}$ denote the norm closures of ${\mathcal A}_0$ 
and ${\mathcal G}_0$, respectively; and
let ${\mathcal{K}}_{X}$ denote the set of all compact operators on $L^2(X)$.
It follows from Theorem 1 in \cite[2.3.4.4]{rempel-schulze} that $\sigma$ and $\gamma$ can be extended to $C^*$-algebra homomorphisms
\[
\bar\sigma:{\mathcal A}\to C(S^*X)\ \ \mbox{and}\ \ \bar\gamma:{\mathcal A}\to C(S^*Y,{\mathcal L}(L^2({\mathbb R}_+))). 
\]
Moreover, by Corollary 2 in \cite[2.3.4.4]{rempel-schulze} and \cite[Theorems 5 and 6]{melo-nest-schrohe1}, we have that:
\begin{equation}
\label{kps}
\ker\bar\gamma\cap\ker\bar\sigma={\mathcal{K}}_{X},\ \ \ \ \ \ker\bar\sigma={\mathcal G},
\end{equation}
and $\bar\gamma$ maps ${\mathcal G}$ onto $C(S^*Y,{\mathcal{K}}_{{\mathbb R}_{+}})$.
In other words, the restriction of the boundary principal symbol to ${\mathcal G}$ gives rise to the 
exact sequence of C$^*$-algebras
%
%
%
%
%
%
\begin{equation}
\label{pses}
0\ \longrightarrow{\mathcal{K}}_{X}\ \longrightarrow \ {\mathcal G}\ {\mathop{\longrightarrow}\limits^{\bar\gamma}}
\ C(S^*Y,{\mathcal{K}}_{{\mathbb R}_{+}})\ \longrightarrow \ 0.
\end{equation}


In Section \ref{sgoi} we use (\ref{pses}) to prove that 
${\mathcal G}$ is isomorphic to the tensor product $\Psi\otimes \mathcal K$
of the $C^*$-closure $\Psi$ of the pseudodifferential operators of order 
zero on $Y$ by the compacts.
For that we need to use trace and Poisson operators.

Similarly as for the singular Green operators, the trace operators and the Poisson operators ($T$ and $K$ in 
(\ref{bdmo})) are, locally, operator-valued-symbol pseudodifferential operators on the variables tangential
to the boundary, given by
\begin{equation}
\label{deft} 
Tu(x^\prime)=(2\pi)^{1-n}\int
e^{ix^\prime\cdot\xi^\prime}t(x^\prime,\xi^\prime,D_n)\hat
u(\xi^\prime)\,d\!\xi^\prime, \ \ u\in C_c^\infty({\mathbb R}^n_+),
\end{equation}
and
\begin{equation}
\label{defk} 
Ku(x^\prime,\cdot)=(2\pi)^{1-n}\int
e^{ix^\prime\cdot\xi^\prime}k(x^\prime,\xi^\prime,D_n)\hat
u(\xi^\prime)\,d\!\xi^\prime,\ \ u\in C_c^\infty({\mathbb R}^{n-1}).
\end{equation}
The mappings $t(x^\prime,\xi^\prime,D_n):L^2({\mathbb R}_{+})\to{\mathbb C}$ and 
$k(x^\prime,\xi^\prime,D_n):{\mathbb C}\to L^2({\mathbb R}_{+})$ are defined, for each 
$(x^\prime,\xi^\prime)\in {\mathbb R}^{n-1}\times{\mathbb R}^{n-1}$, each $v\in L^2({\mathbb R}_{+})$ and
each $\alpha\in{\mathbb C}$, by
\begin{equation}
\label{defst}
t(x^\prime,\xi^\prime,D_n)v=\int\tilde t(x^\prime,y_n,\xi^\prime)v(y_n)\, d\!y_n
\end{equation}
and
\begin{equation}
\label{defsk}
[k(x^\prime,\xi^\prime,D_n)\alpha](x_n)= \alpha \tilde k(x^\prime,x_n,\xi^\prime).
\end{equation}
For each $(x^\prime,\xi^\prime)$, $\tilde t(x^\prime,\cdot,\xi^\prime)$ and $\tilde k(x^\prime,\cdot,\xi^\prime)$
are restrictions to ${\mathbb R}_{+}$ of functions in the Schwartz class on ${\mathbb R}$.
The functions $\tilde t(x^\prime,y_n,\xi^\prime)$ and $\tilde k(x^\prime,x_n,\xi^\prime)$, called the symbol-kernels of 
$T$ and $K$, are smooth and satisfy certain estimates. In the polyhomogenous case, they have asymptotic expansions in 
homogeneous components, whose leading terms we denote by $\tilde t_0$ and $\tilde k_0$, respectively. The
estimates and expansions for $\tilde t$ and $\tilde k$ listed or explained in \cite[Section 1.2]{grubb} are not the 
right ones for our definition of order (and consequent choice of Hilbert
space): we need to shift some of the exponents there by $\pm 1/2$.

The boundary-principal symbols of $T$ and $K$ are
\[
\gamma_T(x^\prime,\xi^\prime)=t_0(x^\prime,\xi^\prime,D_n)\ \ \ \mbox{and}\ \ \   
\gamma_K(x^\prime,\xi^\prime)=k_0(x^\prime,\xi^\prime,D_n), 
\]
defined as in (\ref{defst}) and (\ref{defsk}), except that $\tilde t_0$ and $\tilde k_0$ replace $\tilde t$ and $\tilde k$.
Lastly, the boundary principal symbol of a polyhomogeneous pseudodifferential operator on $Y$ is simply its
usual principal symbol, and we get a *-homomorphism 
\[
\gamma: {\mathcal B}_0\longrightarrow  C^\infty(S^*Y,{\mathcal L}(L^2({\mathbb R}_+)\oplus{\mathbb C})),
\]
where ${\mathcal B}_0$ denotes the set of all polyhomogeneous Boutet de Monvel operators of order and class zero on $X$.
It has a continuous extension to the norm-closure of ${\mathcal B}_0$, but we will not use this fact. 

%
%
%
%

\section{A Product Description of the Singular Green Operators}
\label{sgoi}

\begin{lemma} \label{es1}
There exists a zero-order Poisson operator $K$ such that $K^*K$ is a strictly positive operator on $L^2(Y)$.
\end{lemma}
{\em Proof}:
 It is well-known that the Dirichlet problem 
\[
\left(
\begin{array}{c}\Delta\\\gamma_0\end{array}
\right):
H^2(X)\longrightarrow \begin{array}{c}L^2(X)\\\oplus \\ H^{3/2}(Y)\end{array}
\]
defines a bounded invertible operator.
We denote by $\lambda^{3/2}$ an order reduction of order 3/2 on
$Y$ and by $\Lambda^{-2}$ and order reduction of order $-2$ on $X$.
This gives us an isomorphism
\[
\left(
\begin{array}{c}
\Delta\Lambda^{-2}
\\
\lambda^{3/2}\gamma_0\Lambda^{-2}
\end{array}
\right)
=
\left(
\begin{array}{cc}
1 &0\\0&\lambda^{3/2}
\end{array}
\right)
\left(
\begin{array}{c}\Delta\\\gamma_0\end{array}
\right)
\Lambda^{-2} :L^2(X)\longrightarrow  \begin{array}{c}L^2(X)\\\oplus \\
L^{2}(Y)\end{array}
\]
which is an element of order and class  0  in Boutet de Monvel's calculus.
Its inverse therefore also is in Boutet de Monvel's calculus; it is of the form
$$\left(P_++G\ \ K \right): \begin{array}{c}L^2(X)\\\oplus \\ L^{2}(Y)\end{array}
\longrightarrow L^2(X)$$
with suitable $P, G$, and $K$ of order and class zero.
In particular, $K$ is a right inverse for the trace operator $T=\lambda^{3/2}\gamma_0\Lambda^{-2}:L^2(X)\to L^2(Y)$.
For $v\in L^2(Y)$ we thus have
$$\|v\|_{L^2(Y)}=\|TKv\|_{L^2(Y)}\le \|T\|_{{\mathcal L}(L^2(X), L^2(Y))}\|Kv\|_{L^2(X)}.$$
We then get $\|Kv\|\ge c\|v\|$ for some $c>0$, so that $K^*K$ is strictly
positive.
 \hfill$\Box$

\begin{lemma} \label{es2}
There exist a trace operator of order and class zero $T:L^2(X)\to L^2(Y)$ and a Poisson operator of order zero 
$K:L^2(Y)\to L^2(X)$ such that $TK$ is equal to the identity operator on $L^2(Y)$, $K^*=T$ and $T^*=K$. 
\end{lemma}
{\em Proof}: Let $K_0$ be a zero-order Poisson operator such that 
$K^*_0K_0$ is a strictly positive operator on $L^2(Y)$, and let 
$Q=(K^*_0K_0)^{-1/2}$. $Q$ is a zero-order pseudodifferential operator on $Y$. Take $K=K_0Q$ and $T=QK^*_0$.
\hfill$\Box$

We denote by  $\Psi$  the norm closure of the algebra of all 
polyhomogeneous pseudodifferential operators of order zero on $Y$,
and by $\bar\sigma:\Psi\to C(S^*Y)$ the continuous extension of the 
principal-symbol homomorphism. It is well-known
(this is mentioned in \cite{atiyah-singer} and follows from \cite[Theorem A.4]{kohn-nirenberg}, or from \cite[Theorem 3.3]{hormander66}) that $\bar\sigma$
induces the short exact sequence of $C^*$-algebras
\begin{equation}\label{as}
0\ \longrightarrow{\mathcal{K}}_{Y}\ \longrightarrow \ \Psi\ {\mathop{\longrightarrow}\limits^{\bar\sigma}}
\ C(S^*Y)\ \longrightarrow \ 0,
\end{equation}
where ${\mathcal{K}}_{Y}$ denotes the ideal of compact operators on $L^2(Y)$. 

By Lemma~\ref{es2} a $C^*$-homomorphism $\Xi:\Psi\to{\mathcal G}$ can be defined by 
\[
\Xi(A)=KAT.
\] 
Since $\Xi(A)$ is compact if $A$ is compact, we can use $\Xi$ to couple the
sequences (\ref{pses}) and (\ref{as}). Together they yield the commutative 
diagram of exact sequences of $C^*$-algebras
\begin{equation}
\label{cd}
\def\mapup#1{\Big\uparrow\rlap{$\vcenter{\hbox{$\scriptstyle#1$}}$}}
\begin{array}{ccccccc}
0\longrightarrow&{\mathcal{K}}_X&\longrightarrow& {\mathcal G}  &
{\mathop{\longrightarrow}\limits^{\bar\gamma}}& C(S^*Y,{\mathcal{K}}_{{\mathbb R}_{+}}) &\longrightarrow 0
\\                &\mapup{}&&\mapup{\Xi}&&\mapup{h}&
\\0\longrightarrow& {\mathcal{K}}_{Y} &\longrightarrow&\Psi  &{\mathop{\longrightarrow}\limits^{\bar\sigma}}& C(S^*Y)&\longrightarrow 0
\end{array}.
\end{equation}

\begin{lemma}\label{here}
The homomorphism $\Xi$ imbeds $\Psi$ as a hereditary subalgebra of $\mathcal{G}$. 
\end{lemma}
 {\em Proof}: We have to prove, that if $0 \leq G \leq KAT$ then $G$ is again of the form $KA_1T$ 
with $A_1 \in \Psi$. 
Since $KT$ acts as the identity on $KAT$ 
it also acts as the identity on $G$ and we therefore get $G=KTGKT=K(TGK)T$. \hfill$\Box$ 

 \begin{lemma}\label{full}
 Let
 $$
\def\mapup#1{\Big\uparrow\rlap{$\vcenter{\hbox{$\scriptstyle#1$}}$}}
\begin{array}{ccccccc}
0\longrightarrow& I_1&\longrightarrow& A_1  &
{\mathop{\longrightarrow}\limits^{q_1}} & B_1 &\longrightarrow 0
\\                &\mapup{\phi_1}&&\mapup{\phi_2}&&\mapup{\phi_3}&
\\0\longrightarrow& I_2 &\longrightarrow&A_2  &{\mathop{\longrightarrow}\limits^{q_2}}& B_2&\longrightarrow 0
\end{array}
$$
 be a commutative diagram of short exact sequences, 
where $\phi_1, \phi_2$ and $\phi_3$ are embeddings.
Then $\phi_2$ is full provided  that $\phi_1$ and $\phi_3$ are full. 
 \end{lemma}
 
 {\em Proof}: We have to prove that the two-sided ideal generated by 
$\phi _2(A_2)$ is dense in $A_1$. 
We thus have to prove that to a given $a\in A_1$ and a given $\varepsilon >0$ 
we can find an element $b$ in the twosided ideal generated by the 
image of $\phi_2$ such that $\| a-b \| < \varepsilon$. 
Since $\phi_3$ is full we can find an element $c$ in the twosided ideal generated by $\phi_3(B_2)$ such that $\| q_1(a)-c \| < \frac{\varepsilon}{2}$. The element $c$ can be lifted to an element $b_1$ in the twosided ideal generated by $\phi_2(A_2)$ and we can therefore find an element $d_1\in I_1$ 
such that $\| a-b_1-d_1\| <\frac{\varepsilon}{2}$. 
Since $\phi_1$ is full there is an element $d_2$ in the twosided ideal generated by $\phi_1 ( I_2)$ with $\| d_1-d_2\|  < \frac{\varepsilon}{2}$. As the desired $b$ we can therefore choose $b=b_1+d_2$.\hfill $\Box$

 \begin{theorem}
 The algebras $\mathcal{G}$ and $\Psi \otimes \mathcal{K}$ are isomorphic.
 \end{theorem}
 
 {\em Proof}: By Lemma \ref{here}, the diagram \ref{cd} and Lemma
 \ref{full} the imbedding $\Xi$ is full and hereditary. It follows
 from the remark below Theorem 8 on page 155 in \cite{connes} that
 $\mathcal{G}$ and $\Psi$ are strongly Morita equivalent. By the
 results in \cite{Brown} and \cite{brown-green-rieffel} we have $\mathcal{G}\otimes
 \mathcal{K}$ is isomorphic to $\Psi \otimes \mathcal{K}$. However
 $\mathcal{G}$ is stable since it is the extension of $\mathcal{K}$
 with a stable algebra, namely $C(S^*Y,\mathcal{K}_{\mathbb{R}_+}) $,
(see Proposition 6.12 in \cite{rordam-stable}). This gives the isomorphism.\hfill $\Box$

\section{Pseudodifferential Operators and Groupoids}\label{psi}


Groupoids were introduced in the context of global analysis when
A. Connes showed that in the case of  foliations the index takes
values in a \cstar\ which is defined as the \cstar\ of the holonomy
groupoid of the foliation. He defined a pseudodifferential calculus on
a foliation using the groupoid structure. 

In several papers (\cite{bm-fp, nwx, lmn, bm-jfa}),
generalizations of this approach to a larger class of groupoids were
achieved. One particular aspect of this theory is that, as A. Connes
showed in \cite{connes} for smooth manifolds, it is possible to define the analytic index
using a groupoid, the \textit{tangent groupoid}.

A groupoid is a small category in which all morphisms are
invertible. This means that a groupoid $G$ has a set of units, denoted
by $\g0$, and two maps called \textit{range} and \textit{source}, 
$\xymatrix{G \ar@<.5ex>[r]^{r} \ar@<-.5ex>[r]_{s} &\g0}$.

Two elements $\gamma, \gamma' \in G$ are composable if and only if
$r(\gamma)=s(\gamma')$:
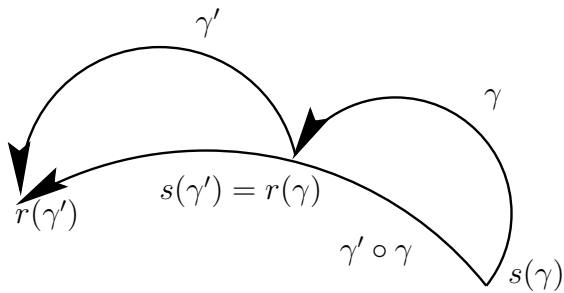
\begin{center}
  \begin{picture}(0,0)%
\epsfig{file=groupd.pstex}%
\end{picture}%
\setlength{\unitlength}{0.000833in}%
\begingroup\makeatletter\ifx\SetFigFont\undefined%
\gdef\SetFigFont#1#2#3#4#5{%
  \reset@font\fontsize{#1}{#2pt}%
  \fontfamily{#3}\fontseries{#4}\fontshape{#5}%
  \selectfont}%
\fi\endgroup%
\begin{picture}(3474,1770)(4576,-3850)
\put(4576,-3436){\makebox(0,0)[lb]{\smash{\SetFigFont{12}{14.4}{\rmdefault}{\mddefault}{\updefault}$r(\gamma')$}}}
\put(5476,-3286){\makebox(0,0)[lb]{\smash{\SetFigFont{12}{14.4}{\rmdefault}{\mddefault}{\updefault}$s(\gamma')=r(\gamma)$}}}
\put(6601,-3661){\makebox(0,0)[lb]{\smash{\SetFigFont{12}{14.4}{\familydefault}{\mddefault}{\updefault}$\gamma' \circ \gamma$}}}
\put(5701,-2236){\makebox(0,0)[lb]{\smash{\SetFigFont{12}{14.4}{\familydefault}{\mddefault}{\updefault}$\gamma'$}}}
\put(7501,-2686){\makebox(0,0)[lb]{\smash{\SetFigFont{12}{14.4}{\familydefault}{\mddefault}{\updefault}$\gamma$}}}
\put(7651,-3811){\makebox(0,0)[lb]{\smash{\SetFigFont{12}{14.4}{\rmdefault}{\mddefault}{\updefault}$s(\gamma)$}}}
\end{picture}

\end{center}


We recall briefly the main aspects of this theory.
Let  $G$ be a Lie  groupoid, which means that it has a smooth
structure. 
Then one can define an  algebra of pseudodifferential operators $\Psi^{\infty}(G)$:
A pseudodifferential operator on $G$ is a $G$-equivariant continuous
family of pseudodifferential operators on the fibers of $G$.

For example, if $M$ is a  manifold without boundary,
    and $G=M\times M$, with set of units $\g0=M$, and range and source
    maps $r(x,y)=x, s(x,y)=y$, and composition $(x,y)(y,z)=(x,z)$,
    then $\Psi^{\infty}(G)$ is the algebra of
    pseudodifferential operators on $M$.

  If $G$ is a Lie group, $\Psi^{\infty}(G)$ is the algebra of
    $G$-equivariant pseudodifferential operators on $G$.

In order to work with singular manifolds, the framework of Lie groupoids
needs to be extended. That was done in \cite{lmn}, where
 the algebras of pseudodifferential operators on continuous
family groupoids, which
are groupoids whose fibers are smooth manifolds, were defined. 

On the algebra of  pseudodifferential operators one can define a
symbol map, $\sigma$. The algebra of order 0 operators can be
completed as a \cstar, denoted by $\overline{\Psi^0}(G)$, and the
symbol map extends to this algebra. The ``regularizing operators'' of
the calculus, which are the operators with trivial symbol, are the elements of the 
\cstar\ of the groupoid, and we have the following Atiyah-Singer
exact sequence:
$$0 \to C^*(G) \to \overline{\Psi^0}(G) \to C(S^*(G)) \to 0,$$
where $S^*(G)$ is the cosphere bundle of the Lie algebroid $A(G)$,
which can be thought of  as a tangent space.





We  next recall in more detail 
the construction of the adiabatic groupoid\ 
$\ad{(Y\times Y)}$ associated with
a smooth manifold $Y$:
$$\ad{(Y\times Y)}=\bigl(TY\times \{0\}\bigl)\ \cup\ \bigl(Y\times Y\times \R^*_+\bigl)$$
with the tangent bundle $TY$ of $Y$.
The groupoid structure is given as follows:
$$r(x,\xi,0)=s(x, \xi, 0)=x,\ (x,\xi,0)\circ (x, \xi',0)=(x, \xi+\xi',0),
$$
$$r(x,y,\lambda)=(x, \lambda),\ s(x,y,\lambda)=(y, \lambda),\ 
(x,y,\lambda)\circ(y,z,\lambda)=(x,z,\lambda),\ \lambda>0.$$

This groupoid is endowed with a differential structure, through an
exponential, in the following way:
\begin{itemize}
\item On $Y\times Y\times \R^*_+$, the structure is that of a product
  of manifolds.
\item Define a map on an open neighborhood $U$ of $TY\times \{0\}$ in
  $TY\times \R_+$, with values in $\ad{(Y\times Y)}$, by
$$
\begin{cases} 
  \psi(x, \xi, \lambda)=(x, \exp_x(-\lambda \xi), \lambda)& \mathrm{\
  if\ }\lambda>0\\
\psi(x, \xi, 0)=(x, \xi,0).&
\end{cases}$$
\end{itemize}
In other terms, the topology is such that a sequence of terms $(x_n,
y_n, \lambda_n)$ of $Y\times Y\times \R^*_+$ converges to $(x,\xi,0)\in TY\times\{0\}$
, if and only if we have locally
$$x_n\to x,\ y_n \to x,\ \frac{x_n-y_n}{\lambda_n} \to \xi.$$
Note that A. Connes' tangent groupoid is just the restriction of  $\ad{(Y\times Y)}$
to $\lambda \in [0,1]$.

The main interest of this groupoid is that it provides a way to define
the analytic index. Consider indeed the decomposition of the groupoid
as an open and a closed subgroupoid, which gives rise to the exact
sequence:
\begin{equation}
  0 \to C^*(Y\times Y\times \R^*_+) \to C^*(\ad{(Y\times Y)}) \to
  C^*(TY) \to 0.
\end{equation}
This simplifies since $C^*(Y\times Y)\simeq \mathcal{K}$, and
$C^*(TY)\simeq C_0(T^*Y)$. A. Connes proved that the boundary map of the 6-terms exact
sequence induced by this extension is nothing but the analytic index 
$$ind_a : K^0_c(T^*Y) \to K_1(C_0(\R)\otimes \mathcal{K})=\Z.$$ 

\section{A Groupoid Associated to the Singular Green Operators}\label{gr-sgo}

Suppose we could identify the $C^*$-closure  of 
Boutet de Monvel's algebra with the 
$C^*$-algebra $\Psi^0(G)$ of pseudodifferential operators on a Lie groupoid $G$.
Then, as pointed out above, the kernel of the principal symbol map
would be isomorphic to $C^*(G)$. 
As the kernel of the principal symbol map in Boutet de Monvel's calculus
consists of the singular Green operators,
we thus wish to identify these with the \cstar\ of a groupoid.

We will actually not
identify them 
with a groupoid \cstar, 
but with an ideal in a groupoid \cstar.

Let us consider the following action of the group $\R^*_+$ on $\ad{(Y\times Y)}$: 
\begin{itemize}
\item On $TY$, $\R^*_+$ acts by dilations: $\lambda . (x,\xi)=(x,
  \lambda\xi)$
\item On $Y\times Y\times \R^*_+$, $\R^*_+$ acts by $\lambda . (x,y,t)=(x,y,\frac{t}{\lambda})$.
\end{itemize}
This is a continuous action: 
If $(x_n, y_n, t_n)$ converges to $(x,\xi)$ (which means that $x_n, y_n
\to x, t_n \to 0, \frac{x_n-y_n}{t_n}\to \xi$), then $\lambda.(x_n,
y_n, t_n)=(x_n, y_n, \frac{t_n}{\lambda}) \to (x, \lambda\xi,0),$ since 
$$\frac{x_n-y_n}{\frac{t_n}{\lambda}}\to \lambda\xi.$$

It is thus possible to construct the semi-direct product 
$G=\ad{(Y\times Y)}\rtimes \R^*_+$ of the adiabatic groupoid by $\R_+$: 
As a set, it  is $\ad{(Y\times Y)} \times
\R^*_+$, with set of units $Y \times \R_+ $, such that:
\begin{itemize}
\item $r(x,y,t,\lambda)=(x,t),\ s(x,y,t,\lambda)=(y,
  \frac{t}{\lambda})$, for $t>0$; 
\item $r(x,\xi,\lambda)=(x,0),
s(x,\xi,\lambda)=(x,0)$, for $t=0$;
\item $(x,y,t,\lambda)(y,z,\frac{t}{\lambda},\mu)=(x,z,t,\lambda\mu);$
\item $(x,\xi,\lambda)(x, \eta, \mu)=(x, \xi+\lambda \eta, \lambda\mu)).$
\end{itemize}
Note that the action of $\R_+^*$ on the adiabatic groupoid 
induces an action on its \cstar\,
and that J. Renault proved in \cite{renault} that
for any locally compact groupoid $\mathcal{G}$ one has
$$C^*(\mathcal{G}\rtimes \R^*_+) \simeq C^*(\mathcal{G}) \rtimes \R^*_+.$$

The evaluation at $t=0$ provides a map $e_0:C^*(G)\to C_0(T^*Y)\rtimes
\R^*_+$. Also, the evaluation at the zero-section $\xi=0$ 
induces a map $r_0:C_0(T^*Y)\rtimes
\R^*_+ \to C(Y)\rtimes \R^*_+$.

But since the action of
$\R^*_+$ on $Y$ is trivial, the latter algebra is just the algebra of
the (regular) product:
$$ C(Y)\rtimes \R^*_+=C_0(Y\times \R^*_+).$$

Let $C=\ker r_0$ and $\mathcal I=\ker r_0 \circ e_0$. 

The kernel of $e_0$ is
$C^*(Y\times Y\times \R^*_+ \rtimes \R^*_+)$. But $\R^*_+ \rtimes
\R^*_+$ is directly isomorphic to the pair groupoid 
$\R^*_+ \times \R^*_+$: To clarify the proof, let us
denote $G_1=\R^*_+ \rtimes \R^*_+$ and $G_2=\R^*_+ \times
\R^*_+$. Then let $\phi:G_1 \to G_2$ be defined by
$$\phi(t,\lambda)=\biggl(t,\frac{t}{\lambda}\biggr).$$
This a morphism of groupoids: 
The composition of  $(t,\lambda)$ with
$\bigl(\frac{t}{\lambda}, \mu\bigr)$ gives $(t, \lambda\mu)$, and 
$$\phi(t, \lambda\mu)=\biggl(t,\frac{t}{\lambda\mu}\biggr)$$
while 
$$\phi(t,\lambda)\circ \phi\biggl(\frac{t}{\lambda},\mu\biggr)=\biggl(t,
\frac{t}{\lambda}\biggr)\circ 
\biggl(\frac{t}{\lambda},\frac{t}{\lambda\mu}\biggr)
=\biggl(t,\frac{t}{\lambda\mu}\biggr).$$

Hence the kernel of $e_0$ is just the algebra of compact operators,
$\mathcal{K}$. 

\medskip
To make this clear, here is the commutative diagram describing this:
\begin{equation}
  \label{diagram}
  \xymatrix{
&&0 \ar[d] & 0 \ar[d]&\\
0 \ar[r] & \mathcal{K} \ar[r] \ar@{=}[d] & \mathcal I  \ar[r] \ar[d] & C \ar[r]
\ar[d]^{j} &0\\
0 \ar[r] & \mathcal{K} \ar[r]& C^*(G)  \ar[r]^<<<<<{e_0} \ar[dr]&
C_0(T^*Y)\rtimes \R^*_+ \ar[r] \ar[d]^{r_0}& 0\\
&&& C_0(Y\times \R^*_+) &\\
}
\end{equation}

We will use this diagram and  extension theory to prove that
 $\mathcal I$ is isomorphic to the algebra of singular Green operators.



\begin{proposition}
  $C$ is isomorphic to $C(S^*Y)\otimes \mathcal{K}$, where $S^*Y$ is the sphere
  bundle in $T^*Y$.
\end{proposition}
\begin{proof}
 First of all notice that $C$ is isomorphic to $C_0(T^*Y\setminus Y)\rtimes
  \R^*_+$: 
Indeed, the exact sequence
$$0 \to C_0(T^*Y\setminus Y) \to C_0(T^*Y) \to C(Y) \to 0$$
induces the exact sequence
$$0 \to C_0(T^*Y\setminus Y)\rtimes \R^*_+ \to C_0(T^*Y)\rtimes \R^*_+ \to
  C(Y)\rtimes \R^*_+ =C_0(Y\times \R^*_+)\to 0.$$

But $T^*Y\setminus Y \simeq S^*Y\times \R^*_+$, so that $C_0(T^*Y\setminus
Y)\simeq C(S^*Y) \otimes C_0(\R^*_+)$ and 
$$C_0(T^*Y\setminus Y)\rtimes
\R^*_+ \simeq C(S^*Y) \otimes C_0(\R^*_+)\rtimes
\R^*_+.$$

Now  
$$C_0(\R^*_+)\rtimes \R^*_+ 
\simeq C^*(\R^*_+)\rtimes \R^*_+ 
\simeq C^*(\R^*_+\rtimes \R^*_+)
\simeq C^*(\R^*_+\times \R^*_+)
\simeq \mathcal{K}$$
where we again used Renault's result for the second isomorphism 
and  the isomorphism of $\R^*_+\rtimes \R^*_+$ 
with the pair groupoid $\R^*_+\times \R^*_+$
for the third.
This ends the proof.

\end{proof}

\section{Identification of the Ideal  with the 
Singular Green Operators}\label{main}

We have just shown that $\mathcal I$ is an extension of $C_0(S^*Y)\otimes \mathcal{K}$ by
$\mathcal{K}$, and this is also the case for the algebra of \sgo s. The
main result is the following:

\begin{theorem}\label{isom}
  The \cstar\ $\mathcal I$ is isomorphic to $\Psi\otimes \mathcal{K}$.
\end{theorem}
\begin{proof}
 
We shall prove that the extensions
$$0 \to \mathcal{K} \to \mathcal I \to C(S^*Y)\otimes \mathcal{K} \to 0$$
 and 
$$0 \to \mathcal{K} \to \Psi \otimes \mathcal{K} \to C(S^*Y)\otimes \mathcal{K} \to 0$$
satisfy the conditions of a theorem by Voiculescu which we recall now.

D. Voiculescu proved in \cite{voiculescu-extensions} (look also at the
survey \cite{skandalis-expositiones}, Theorem 10.9) that if
two extensions $0 \to \mathcal{K} \to D_1 \to A\to 0$ and  $0 \to
\mathcal{K} \to D_2 \to A \to 0$, are such that:
\begin{itemize}
\item $D_1$ and $D_2$ define the same element in $Ext(A)$,
\item $D_1$ and $D_2$ are not unital,
\item $\mathcal{K}$ is essential in $D_1$ and in $D_2$,
\end{itemize}
then $D_1$ and $D_2$ are isomorphic. Recall that an ideal $J$ of $D$ is essential
if and only if for every $x\in D$, $x\neq 0 \Rightarrow  \exists y \in J, xy\neq 0$.

Let us apply  this important result  in our context.

First of all, $\mathcal I$ and $\Psi\otimes \mathcal{K}$ are
non-unital since their quotients by the compacts are isomorphic to
$C(S^*Y)\otimes \mathcal{K}$, which is non-unital.

 Since the algebra
of compact operators on a Hilbert space $\mathcal{H}$ is essential in
any \cstar\ included in $\mathcal{L}(\mathcal{H})$, we obtain that
$\mathcal{K}_X$ is essential in $\Psi \otimes \mathcal{K}_{\R_+}$.

For the other algebra, notice that $C^*(Y\times Y \times \R^*_+)$ is
an essential ideal of $C^*(\ad{(Y\times Y)})$, thus its crossed product by
$\R^*_+$ is also an essential ideal of  $C^*(\ad{(Y\times Y)})\rtimes \R^*_+$,
hence of $\mathcal I$ (see \cite{kusuda}).

It remains to show that the extensions give rise to the same element
of $Ext(C(S^*Y)\otimes \mathcal{K})$. But since $C(S^*Y)\otimes
\mathcal{K}$ is separable, $KK_1(C(S^*Y)\otimes \mathcal{K}, \C)$ is
isomorphic to the group of invertibles of $Ext(C(S^*Y)\otimes
\mathcal{K})$, thanks to a result of Kasparov
(\cite{kasparov-izv}). The \cstar\ $C(S^*Y)\otimes \mathcal{K}$ being
nuclear, $Ext(C(S^*Y)\otimes
\mathcal{K})$ is actually a group, thus it is isomorphic to
$KK_1(C(S^*Y)\otimes \mathcal{K}, \C)$.

The element  $i_S \in KK_1(C(S^*Y)\otimes \mathcal{K}, \C)\simeq Ext(C(S^*Y)\otimes \mathcal{K})$ associated to the
extension
$$0 \to \mathcal{K}_X \to \Psi \otimes \mathcal{K}_{\R_+} \to C(S^*Y)\otimes \mathcal{K}_{\R_+} \to 0$$
provides a map $K_1(C(S^*Y)\otimes \mathcal{K}_{\R_+})\simeq
K_1(C(S^*Y)) \to K_0(\C)=\Z$, which is the analytic index. 

For the class of the extension
$$0 \to \mathcal{K} \to \mathcal I \to C(S^*Y)\otimes \mathcal{K} \to 0,$$
let us consider first the extension
$$0 \to C^*(Y\times Y \times \R^*_+) \to C^*(\ad{(Y\times Y)}) \to
C_0(T^*Y) \to 0,$$ 
whose class is denoted by $i_T \in KK_1(C_0(T^*Y),
\mathcal{K} \otimes C_0(\R^*_+))\simeq KK_0(C_0(T^*Y), \mathcal{K})$.


It induces 
the extension
 $$0 \to C^*(Y\!\times\! Y \!\times\! \R^*_+)\rtimes \R^*_+ \simeq \mathcal{K} \to
 C^*(\ad{(Y\!\times\! Y)})\rtimes \R^*_+ \simeq C^*(G) \to
 C_0(T^*Y)\rtimes \R^*_+ \to 0$$
whose class %
is denoted by $\alpha \in KK_1(C_0(T^*Y)\rtimes \R^*_+,\mathcal{K})$. 

The relation between
$i_T$ and $i_S$
is made clear by considering the
following  exact sequence
$$0 \to C_0(T^*Y) \to C(B^*Y) \to C(S^*Y) \to 0$$
where $B^*Y$ is the ball bundle over $Y$. 
Its class is denoted by $\psi \in KK_1(C(S^*Y), C_0(T^*Y))$, and one
has the well-known equality:
$$i_S=\psi i_T.$$

For the convenience of the reader, we now recall the diagram (\ref{diagram}):
$$  \xymatrix{
0 \ar[r] & \mathcal{K} \ar[r] \ar@{=}[d] & \mathcal I  \ar[r] \ar[d]
& C \ar[r]
\ar[d]^{j} &0  & (\partial)\\
0 \ar[r] & \mathcal{K} \ar[r]& C^*(G)  \ar[r]^<<<<<{e_0} \ar[dr]&
C_0(T^*Y)\rtimes \R^*_+ \ar[r] \ar[d]^{r_0}& 0 & (\alpha)\\
&&& C_0(Y\times \R^*_+) & &\\
}$$
Let us denote the class of the
first sequence by $\partial \in KK_1(C, \mathcal{K})$; it is thus given by the Kasparov product:
$$ \partial= j_* \alpha.$$

Let us make the relation between $j_*$ and $\varphi$ precise.

Consider  the following commutative diagram: 
$$
    \xymatrix{
0 \ar[r] & C_0(T^*Y\setminus Y) \ar[r] \ar[d] &
C_0(B^*Y\setminus Y) \ar[r] \ar[d] & C(S^*Y) \ar[r] 
\ar@{=}[d] &0\\ \label{eq1}
0 \ar[r] & C_0(T^*Y) \ar[r] \ar[d]& C(B^*Y) \ar[d] \ar[r] &
C(S^*Y) \ar[r] & 0\\
&C(Y) \ar@{=}[r] &C(Y)& &\\ \label{eq2} 
}
$$
The first exact sequence actually decomposes as
$$0 \to C(S^*Y)\otimes C_0(\R^*_+) \to  C(S^*Y)\otimes C_0(\R^*_+ \cup
\{\infty\}) \to  C(S^*Y) \to 0$$
so that its $KK_1$-class is the identity of $KK_1(C(S^*Y),  C(S^*Y)\otimes C_0(\R^*_+))$.

There is an action of $\R^*_+$ on each algebra of the previous diagram, which is trivial on
$C(S^*Y)$ and $C(Y)$. This gives the following:
$$
    \xymatrix{
0 \ar[r] & C_0(T^*Y\!\setminus\! Y)\!\rtimes \R^*_+ \ar[r] \ar[d]^{j} &
C_0(B^*Y\!\setminus\! Y)\!\rtimes \R^*_+ \ar[r] \ar[d] & C(S^*Y\times \R^*_+) \ar[r] 
\ar@{=}[d] &0 &(\partial_1)\\ \label{eq1a}
0 \ar[r] & C_0(T^*Y)\!\rtimes \R^*_+ \ar[r] \ar[d]& C(B^*Y)\rtimes \R^*_+ \ar[d] \ar[r] &
C(S^*Y\times\R^*_+) \ar[r] & 0 &(\varphi)\\
&C(Y\times\R^*_+) \ar@{=}[r] &C(Y\times\R^*_+)& & &\\ \label{eq2a} 
}
$$

Denote by  $\partial_1$ (resp. $\varphi$) the class of the
first (resp.  second)  exact sequence of this diagram, and by $j_* \in
KK(C_0(T^*Y\setminus Y), C_0(T^*Y))$ the element
induced by $C_0(T^*Y\setminus Y) \to C_0(T^*Y)$. One has thus the equality
$$\varphi=  \partial_1 j_*,$$

so that
$$\partial=j_* \alpha= \partial_1^{-1} \varphi \alpha.$$

But $ \varphi \alpha$ is the image of $\psi i_T$ under the Thom-Connes
isomorphism, and $\partial_1$ is also a Thom-Connes element in $KK$-theory. 
Hence the classes in $KK_1$ of the extensions of
$\mathcal{I}$ and of $\Psi \otimes \mathcal{K}$ are the
same. Voiculescu's theorem implies that these algebras are isomorphic.
\end{proof}

\begin{corollary} 
   The algebra of singular Green operators is isomorphic to $\mathcal I$, an
  ideal of $C^*(\ad{(Y\times Y)}\rtimes \R^*_+)$.
\end{corollary}
\begin{proof}
  This is a direct consequence of Theorem \ref{isom}, since the
  algebra of singular Green operators is isomorphic to $\Psi \otimes \mathcal{K}$.
\end{proof}


\end{document}